\documentclass[11pt]{amsart}
\usepackage{amsmath,amsfonts,amsthm, amscd, mathrsfs, amssymb, txfonts}
\usepackage{xy}
\numberwithin{equation}{section}
\newtheorem{theorem}{Theorem}[section]
\theoremstyle{definition}

\newtheorem{lemma}[theorem]{Lemma}
\newtheorem{example}[theorem]{Example}
\newtheorem{proposition}[theorem]{Proposition}

\newtheorem{conjecture}[theorem]{Conjecture}
\theoremstyle{definition}
\newtheorem{remark}[theorem]{Remark}

\newcommand{\Ker}{\textrm{Ker}}

\newcommand{\Aut}{\textrm{Aut}}

\newcommand{\Id}{\textrm{Id}}

\newcommand{\Fix}{\textrm{Fix}}
\newcommand{\Tr}{\textrm{Tr}}
\newcommand{\Alb}{\textrm{Alb}}
\begin{document}

\title{Fixed Points of Smooth Varieties with Kodaira Dimension Zero}
\author{Adam T. Ringler}
\date{August 27, 2007}
\address{Department of Mathematics and Statistics, University of New Mexico, Albuquerque, New Mexico 87131}
\email{ringler@unm.edu} 
\maketitle

\begin{abstract}
In this paper, we study the growth of the number of fixed points from iterating an endomorphism of an abelian variety.  Using the estimates obtained on an abelian variety, we are able to extend the results to endomorphisms on varieties of Kodaira dimension zero and more generally their periodic subvarieties. 
\end{abstract}

\section{Introduction}
\label{intro}
A fundamental problem in the theory of dynamical systems is to study the fixed points of a self map $\phi : X \to X$ in various settings.  For example, when $X$ is a smooth manifold and $\phi$ is a smooth map, Shub and Sullivan \cite{Shub} raised, roughly stated, the following question: Does
$$
\# \Fix (\phi^l) = \# \{ P \in X : \phi^l (P) =P \}
$$
have exponential growth as a function of $l$?  Turning to a more algebraro-geometric setting we can ask the same question when $X$ is a variety over an algebraically closed field $k$ of characteristic $0$ and $\phi : X \to X$ is an endomorphism.  In this setting Zhang \cite[Conjecture 2.1.1]{Zhang} has conjectured the following:

\begin{conjecture}
Let $X$ be a projective variety over $\mathbb{C}$ and $\varphi : X \to X$ an endomorphism with $\varphi^*L \cong L^{\otimes q}$ for some ample line bundle $L$ on $X$.  Suppose that $Y \subset X$ is a subvariety of dimension $r$ with $\varphi^m (Y) = Y$.  Then
$$
\# \{ P \in Y(\mathbb{C}) : \varphi^{ml} (P)=P \} = q^{rml} (1+ o(1))
$$
as $l \to \infty$.
\end{conjecture}

Motivated by the above problems, we now turn to the results of this paper.  We first look at the case when $A$ is an abelian variety and $\phi : A \to A$ is an endomorphism satisfying a mild polarization hypothesis.  In this setting we are able to prove that $\# \Fix (\phi^l)$ grows exponentially as a function of $l$.  

More preciesely, let $A_i$ be a simple abelian subvariety of $A$ of dimension $g_i$ for $i=1,\dots, m$.  Suppose that for each $A_i$ there exists an ample line bundle $L_i$ on $A_i$ such that $\phi^* L_i \cong L_i ^{\otimes q_i}$ for some $q_i >1$ and $\phi (A_i) = A_i$.  In this case we are able to prove:

\begin{theorem}\label{a}
Let $A$, $\phi$, and $A_1 ,\dots , A_m$ be as above.  Suppose that $A$ is isogenous to $A_1 \times \cdots A_m$.  Then
$$
\Fix (\phi^l) =\prod_{i=1} ^m (q_i ^l -1)^{g_i} +C
$$
where $C$ is a constant that depends only on $A$ and the product is taken over $A_1, \dots, A_m$ simple abelian subvarieties.
\end{theorem}

It is interesting to note that we do not require that ample line bundle $L_i$ on $A_i$  to be the restriction of an ample line bundle on $A$.  By dropping the global polarization hypothesis, we get a more refined estimate for a larger class of maps, than we would had we assume $\phi^* L \cong L^{\otimes q}$ for some ample line bundle $L$ on $A$. 

As an application of our fixed point estimates on abelian varieties, we are able to make a small step towards the above conjecture.

\begin{theorem}\label{b}
Let $A$ be an abelian variety over an algebraically closed field $k$ of characteristic $0$.  Let $\phi : A \to A$ be a morphism and $L$ an ample line bundle on $A$ with $\phi^* (L) \cong L^{\otimes q}$ for some $q >1$.  If $Y \subset A$ is a smooth subvariety of dimension $r$ with $\phi^m(Y) =Y$ for some $m >0$,  then
$$
\# \{ P \in Y(k) : \phi^{ml}(P)=P \} = q^{rml}(1+o(1))
$$
as $l \to \infty$.
\end{theorem}

Exploiting the above fixed point counts for abelian varieties one is able to extend the previous theorem, at the cost of a constant, to the case where $\phi : X \to X$ and $X$ is a smooth projective variety of Kodaira dimension $0$. 

\begin{theorem}\label{c}
Let $X$ be a smooth projective variety over $\mathbb{C}$ of Kodaira dimension zero.  Let $L$ be an ample line bundle on $X$ and $\phi : X \to X$ an endomorphism with $\phi^* L \cong L^{\otimes q}$ for some $q >1$.  If $Y \subset X$ is a smooth subvariety of dimension $r$ with $\phi^{m} (Y) = Y$ for some $m >0$, then there exists a constant $C$ which depends only on $Y$ and $X$ such that
$$
\# \{ P \in Y(\mathbb{C}) : \phi^{ml} (P) = P \} \geq Cq^{rml} (1+o(1))
$$
as $l \to \infty$.
\end{theorem}

The main difficulty in the above estimates arises from the fact that the fixed points are not counted with multiplicity.  This makes it difficult to apply the Lefschetz fixed point theorem directly or to use B\'ezout type arguments.  

As was pointed out by Zhang, if $\phi = [m]$ is the multiplication by $m$ map, then theorem \ref{b} reduces to a simple result about torsion points on an abelian variety.  For sake of completeness we will give the argument, in this elementary case, in section 1.

The outline of the paper is as follows: in section $1$ we recall some background results which will be used later on and give the argument for the proof of our first theorem when $\phi : A \to A$ is the multiplication by $m$ map.  In section $2$ we prove theorem \ref{a}.  In section $3$  we prove theorem \ref{b} by reducing to a case where we may apply theorem \ref{a}.  Finally, in section $4$ we apply the conclusion of theorem \ref{b} to prove theorem \ref{c}.


\section{Background and Conventions}\label{sec:1}
We now fix some notation and conventions.  Let $k$ be an algebraically closed field of characteristic $0$.  We will assume all varieties are integral and morphisms are defined over $k$.  In section $4$ we will restrict to the case $k = \mathbb{C}$.  

Let $A$ be an abelian variety of dimension $g$.  If $\phi : A \to A$ is an isogeny, then its degree is the cardinality of its kernel.  Recall that to every isogeny $\phi$ there exists a unique dual isogeny $\hat{\phi}$ such that $\phi \circ \hat{\phi} = \hat{\phi} \circ \phi =[m]$ for some $m$, where $[m] : A \to A$ denotes the multiplication by $m$ map, which is the map that sends $P \mapsto mP$ under the group law on $A$.  If $L$ is an ample symmetric line bundle on $A$ that is, $L \cong [-1]^* L$, then the pull-back of $L$ by $[m]$ satisfies $[m]^* L \cong L^{\otimes m^2}$.  We refer the reader to \cite{Milne1} for other basic results on abelian varieties.  

We will need a standard result about the intersection number of divisors upon pull-back, see \cite[Example 2.4.3]{Fulton}.  Namely, let $\varphi : X \to Y$ be a finite morphism and $X$ and $Y$ smooth projective varieties over $k$ of dimension $n$.  Let $D_1, \dots, D_n$ be divisors on $Y$.  Then
\begin{equation}\label{inter}
(\varphi^* (D_1) \cdots  \varphi^* (D_n)) = \deg (\varphi) (D_1  \cdots  D_n),
\end{equation}
where $\cdot$ denotes the intersection. 

Restricting to the case where $Y=X$, so that $\varphi: X \to X$ is an endomorphism, we denote by $\Fix (\varphi)$ the set of fixed points counted without multiplicity.  In this case, see theorem \ref{fakh} below, the Kodaira dimension of $X$ can not be positive.  In other words, if $K_X$ denotes the canonical divisor, then the linear system $|mK_X|$ does not give rise to any map from $X$ to a projective space with positive dimensional image for any $m>0$.  

Following the definition given in \cite{Debarre}, we say a divisor $D$ on $X$ is \emph{generically nef} if $D \cdot H_1\cdots H_{n-1} \geq 0$ for all ample divisors $H_1, \dots, H_{n-1}$ on $X$.  This will be used to apply the following lemma \cite[Lemma 3.18]{Debarre}

\begin{lemma}\label{deblem}
Let $Y$ be a normal variety, $A$ an abelian variety, and $\psi : Y \to A$ a finite map.  If $-K_Y$ is generically nef, then $Y$ is an abelian variety.
\end{lemma}

We will also need the following elementary lemma which will be useful on two occasions.  

\begin{lemma}\label{pointwise}
Let $X$ be smooth projective variety and $\phi : X \to X$ an endomorphism such that $\phi^* L \cong L^{\otimes q}$ for some $q >1$.  Then $\phi$ has finitely many fixed points.
\end{lemma}

\begin{proof}
It suffices to check that $\phi$ does not pointwise fix a positive dimension subvariety.  Suppose this is not the case so that $\phi(V) = V$ pointwise, where $V \subset X$ is a positive dimensional subvariety.  Then $\phi_* (V) =q^{\dim (V)}[V]$, where $[V]$ denotes the cycle associated to $V$.  However, $\phi|_V : V \to V$ is the identity so $[V]=q[V]$, which is impossible since $q >1$.
\end{proof}

We recall a theorem of Fakhruddin, which gives strong restrictions on the geometry of a polarized endomorphism $\phi : X \to X$, when $X$ has non-negative Kodaira dimension, denoted $\kappa(X) \geq 0$.

\begin{theorem}\label{fakh}
Let $X$ be a smooth projective variety with $\kappa(X) \geq 0$, $L$ an ample line bundle, and $\phi : X \to X$ an endomorphism with $\phi^* L \cong L^{\otimes q}$ for some $q >1$.  Then $m K_X \sim 0$ for some $m  \geq 1$ and $\phi$ is \'etale.  Furthermore, $X$ is isomorphic to the quotient of an abelian variety $A$ by a finite group $G$ acting fixed point freely on $A$ and $X$ contains no positive dimensional subvarieties with finite algebraic fundamental group.
\end{theorem}

\begin{proof}
See lemma 4.1 and theorem 4.2 of \cite{Fakhruddin}.
\end{proof}

By \cite[Expos\'e XI, Corollarie 1.3]{Grothendieck}, a variety $X$ satisfying the conditions  of the previous theorem can not contain any rational curves, because otherwise $X$ would contain a subvariety with a finite algebraic fundamental group. 

Using the previous theorem Zhang has given a uniruledness criterion in all dimensions for an endomorphism of a smooth projective variety with negative Kodaira dimension.

\begin{proposition}\label{unir}
Let $X$, $\phi$, and $L$ be as in the previous theorem and suppose further that $\kappa(X) <0$.  Then $X$ is uniruled.  
\end{proposition}

\begin{proof}
See \cite[Proposition 2.2.1]{Zhang}.
\end{proof}

As stated in the introduction, we give an argument for the case where $\phi : A \to A$ is the multiplication by $m$ map on an abelian variety $A$ of dimension $g$.  This case has already been observed by Zhang in \cite{Zhang}.  However, as theorem \ref{b} follows the basic outline of this elementary case we give the argument here.

In the case where $\phi= [ m]$, we see that if $P \in A(k)$ is a fixed point of $\phi^l$ for $l >0$, then $P$ is in  the $m^l -1$ torsion of $A$, denoted $A[m^l-1]$.  By a standard result on abelian varieties, see \cite[Theorem A.7.2.7]{Silverman}, we know that
$$
\# A[m^l -1] = (m^l-1)^{2g}.
$$
In this case, as noted by Zhang, periodic subvarieties will be translates of abelian subvarieties by torsion points.  For a proof of this fact, see \cite[Lemma 4.1]{Faltings}.  Therefore, it is easy to conclude theorem \ref{b} in this case.  

The main complication in our more general setting is that we no longer have an explicit count for the fixed points of $\phi$ with respect to some torsion subgroup of $A$.  In our case, the set of fixed points could be a proper subset of some torsion subgroup.  One also needs to give a more geometric proof that the only possible periodic subvarieties will be translates of abelian subvarieties.

We conclude this section by giving an example of an isogeny, which is not the multiplication by $m$ map for which our results hold.  As theorem \ref{a} only requires that certain simple abelian varieties be polarized, it is also possible to give more elaborate examples.

\begin{example}
The following example comes from \cite[Exercise A.7.9]{Silverman}.  Let $\phi : A \times A \to A \times A$ be given by $(x,y) \mapsto (x+y, x-y)$ and $D$ an ample symmetric line bundle on $A$.  Then $\phi$ has degree $2^{g}$ and 
$$
\phi^*(D \times A + A \times D) \sim 2(D \times A + A \times D)
$$
so the conditions of our main result are satisfied.
\end{example}

\section{Periodic Points of $A$}\label{sec:2}
In this section we prove theorem \ref{a}.  The result will follow from three lemmas.  First we prove that if $\phi$ is an isogeny, then so is $\phi -\Id_A$.  Then as $\Ker(\phi-\Id_A)$ is the set of fixed points of $\phi$ and is the same as $\deg(\phi-\Id_A)$, it suffices to compute the degree of $\phi- \Id_A$ on the simple abelian subvarieties of $A$, which is the content of lemma \ref{degprod}.

\begin{lemma}
Let $A$ be an abelian variety of dimension $g$.  If $\phi$ is an isogeny on $A$ such that $\phi^* L \cong L^{\otimes q}$ for some ample line bundle $L$ on $A$ and $q >1$, then $\phi - \Id_A$ is also an isogeny on $A$.
\end{lemma}


\begin{proof}
Clearly, $\phi- \Id_A$ is a homomorphism of $A$ and since $\phi - \Id_A$ is an endomorphism, it suffices to show that $\Ker(\phi - \Id_A)$ is finite.  Suppose that $\Ker(\phi-\Id_A)$ is not finite, then $\phi$ must pointwise fix a positive dimensional subvariety, which would contradict lemma \ref{pointwise}.
\end{proof}

We now prove a basic result on the intersection numbers of divisors in products of a variety.

\begin{lemma}\label{proddiv}
Let $X^r$ be the product of $r$ copies of a smooth projective variety $X$ of dimension $n$.  Let $D_i$ be a divisor on the $i$-th copy of $X$, denoted $X_i$.  Let $F_i =(X^{i-1} \times D_i \times X^{r-i})$, which is a divisor on $X^r$.  Then
$$
(F_1 + \dots + F_r)^{rn} = ( r! D_1 \times \cdots \times  D_r)^{n},
$$
where $( \cdot)^m$ denotes the self-intersection $m$-times.
\end{lemma}
\begin{proof}
First recall $(D_i)^{n+1} =0$, since each $X_i$ has dimension $n$.  Since $D_1 \times \cdots \times D_r$, will be the only term where an intersection of the form $(D_i)^{n+1}$ does not occur, it suffices to compute the coefficient.  Since 
$$
(F_1 \cdot F_2 \cdots F_r)^n =(D_1 \times \cdots \times D_r)^n
$$
occurs $r!^n$ times in the expansion of $(F_1 + \dots +F_r)^{rn}$, our lemma follows.
\end{proof}

\begin{lemma}\label{degprod}
Let $A^r$ be the product of $r$-copies of a simple abelian variety $A$ of dimension $g$.  Let $\phi : A^r \to A^r$ be an isogeny  with $\phi|_{A_i} ^* L_i \cong L_i ^{\otimes q_i}$, where $A_i$ denotes the $i$-th copy of $A$, $L_i$ is an ample line bundle on $A_i$, and $q_i >1$ for $i=1, \dots, r$.  Then
$$
\deg(\phi - \Id_{A^r}) =\prod_{i=1} ^r (q_i-1)^g.
$$
\end{lemma}

\begin{proof}
Let $\psi : A^r \to A^r$ be the dual isogeny of $\phi - \Id_{A^r}$.  Then
\begin{equation}\label{first}
\psi \circ ( \phi - \Id_{A^r}) = [m]
\end{equation}
for some $m >0$ which we may rewrite as
$$
\psi = \psi \circ \phi - [m].
$$
Putting this back into (\ref{first}) yields
$$
\psi \circ (\phi - \Id_{A^r}) =\psi \circ \phi^2 - \psi \circ \phi - [m] \circ \phi +[m].
$$
Therefore, we obtain
$$
\psi \circ \phi^2 - \psi \circ \phi - [m] \circ \phi = [0].
$$
As $(\phi - \Id_A) \circ \psi = \psi \circ (\phi - \Id_A)$, we see that $\phi \circ \psi = \psi \circ \phi$.  Let $D_i$ be a divisor on $A_i$ that satisfies $\mathcal{L}(D_i)  \cong L_i$.  Recalling the notation from the previous lemma, set
$$
F_i = A^{i-1} \times D_i \times A^{r-i}  \textrm{ \    \ for \    \ } i =1, \dots, r. 
$$
Then $F_i$ is a divisor on $A^r$.  If one of the $F_i$ is not symmetric, then replace $D_i$ by $D_i+ [-1]^* D_i$, where it appears in each of the divisors $F_i$, symmetric or not.  Then each $F_i$ will be symmetric as $[-1]$ is an isomorphism.  Furthermore, each $F_i$ will satisfy $\phi^* F_i  \sim q_i F_i$ as $\phi$ commutes with $[-1]$.  So we may assume that each $F_i$ is symmetric for $i =1, \dots r$.  We then have
$$
(\phi^2 \circ \psi)^* F_i - ( \psi \circ \phi)^* F_i - ([m] \circ \phi)^* F_i \sim 0.
$$
As $\phi$ and $\psi$, commute we obtain
$$
\psi^* (q_i ^2 F_i) - \psi^* (q_i F_i) - [m]^* (q_i D_i) \sim 0.
$$
Since $F_i$ is symmetric ($F_i \sim [-1]^*F_i$), then as $[m]^*F_i \sim m^2  F_i$, we obtain
\begin{equation}\label{second}
\psi^* ( q_i(q_i-1)F_i) - m^2 q_i F_i \sim 0.
\end{equation}
 Using (\ref{inter}), we obtain
$$
(\psi^* ((q_1^2 -q_1)F_1 + \dots + (q_r ^2-q_r)F_r))^{rg} = \deg(\psi)((q_1^2-q_1)F_1 + \dots +(q_r^2 -q_r)F_r)^{rg}
$$
where $(\cdot)^{rg}$ denotes the self-intersection of the divisor in question $rg$-times.  Then applying lemma \ref{proddiv}  and (\ref{second}), we obtain
$$
(m^{2r} q_1 \cdots q_r)^{g} (D_1 \times \cdots  \times D_r)^g= \deg(\psi)((q_1^2 -q_1) \cdots (q_r ^2 -q_r))^{g}( D_1 \times \cdots \times D_r )^g.
$$
As each $D_i$ is ample for $i=1 ,\dots , r $ we have
$$
 (D_1 \times \cdots \times D_r) ^g >0
$$
which implies that
$$
m^{2rg} = \deg(\psi) ((q_1-1) \cdots (q_r -1))^g.
$$
The lemma now follows by recalling that $\deg(\psi)\deg(\phi- \Id_{A^r}) = m^{2rg}$ as $\psi$ is the dual isogeny of $\phi-\Id_{A^r}$.
\end{proof}

Using the previous lemma, we are now ready to prove theorem \ref{a}. \newline


\noindent{\it Proof of theorem \ref{a} }  First we recall some notation.  Let $\phi: A \to A$ be an endomorphism of an abelian variety $A$.  let $A_i$ be an abelian subvariety of $A$ of dimension $g_i$ for $i=1,\dots, m$.  Suppose that for each $A_i$ there exists an ample line bundle $L_i$ on $A_i$ such that $\phi^* L_i \cong L^{\otimes q_i}$ for some $q_i >1$.  We further suppose that $A$ is isogenous to $A_1 \times  \cdots \times A_m$.  For the purpose of counting fixed points, up to a constant $C$, it suffices to assume, possibly after reindexing, that
$$
A = A_1 ^{r_1} \times \cdots \times A_m ^{r_m},
$$
where $A_1, \dots, A_m$ are distinct simple abelian subvarieties and $A_i$ has dimension $g_i$.  Then we may write $\phi =  \tau_a \circ (\phi_1 \times \cdots \times \phi_m)$, where each $\phi_i : A_i ^{r_i} \to A_i ^{r_i}$ is an isogeny for $i=1, \dots, m$ and $\tau_a : A \to A$ is the translation by $a \in A(k)$.  As $\phi_1 \times \cdots \times \phi_m$ is unramified and so is $\phi - \Id_A$, we see that
$$
\# (\phi_1 \times \cdots \times \phi_m - \Id_A)^{-1} (-a) = \# ( \phi_1 \times \cdots \times \phi_m -\Id_A)^{-1} (0).
$$
Therefore, we see that the number of fixed points of $\phi$ is 
$$
 \# \Fix (\phi) = \deg (\phi_1 \times \cdots \times \phi_r - \Id_A).
$$
Hence,
$$
\# \Fix (\phi) = \# \Fix(\phi_1) \cdots \# \Fix (\phi_r) +C,
$$
where $C$ is the previously introduced constant and independent of $\phi$.  Let $A_{ij}$ denote the $j$-th copy of $A_i$ in $A_i ^{r_i}$.  Then there exists an ample line bundle $L_{ij}$ on $A_{ij}$ such that $\phi|_{A_{ij}} ^* L_{ij} \cong L_{ij} ^{\otimes q_{ij}}$.  Therefore, it suffices to show
$$
\# \Fix (\phi_i ^l) = \prod_{j=1} ^{r_m} (q_{ij} ^l -1)^{g_i}.
$$
However, this follows from lemma \ref{degprod} by replacing $\phi_i$ with $\phi_i ^l$ since 
$$
\# \Fix (\phi_i ^l) = \deg (\phi_i ^l - \Id_{A_i ^{r_i}})
$$
which completes the proof.

It is interesting to note that the above theorem shows where the fixed points are coming from.  Namely, if $A_i$ is one of the simple abelian subvarieties in the statement of the theorem, then it contributes approximately $(q_i^l -1)^{g_i}$ fixed points.  Then taking the product of all the simple abelian subvarieties gives the total number of fixed points up to a constant.

\section{Periodic Points of Subvarieties of Abelian Varieties}\label{sec:3}

We now prove theorem \ref{b}.  The idea is to use theorem \ref{a} to give a uniform asymptotic estimate when the endomorphism $\phi : A \to A$ satisfies $\phi^* L \cong L^{\otimes q}$ for some ample line bundle $L$ on $A$, where $q >1$.  We then study the possible smooth subvarieties $Y \subset A$ such that $\phi^m (Y) = Y$ for some $m >0$.  As one would expect, it turns out that $Y$ must be a translate of an abelian subvariety of $A$, which is the content of lemma \ref{perab}.  

\begin{proposition}\label{third}
Let $A$ be an abelian variety of dimension $g$ and $L$ an ample line bundle on $A$.  Let $\phi : A \to A$ be an endomorphism with $\phi^* (L) \cong L^{\otimes q}$ for some $q >1$.  Then
$$
\# \{ P \in A(k) : \phi^l (P) =P \} =q^{gl} (1+o(1))
$$
as $l \to \infty$.
\end{proposition}

\begin{proof}
First notice that it suffices to prove the claim under the assumption that $\phi$ is an isogeny.  Recall that any endomorphism of an abelian variety can be written as the composition of a translation $\tau_a$ for some $a \in A(k)$ and an isogeny $\psi$ so that $\tau_a \circ \psi = \phi$.  If $P$ is a fixed point of $\phi$, then $\psi (P) +a =P$.  By theorem \ref{fakh}, $\phi - \Id_A$ is unramified so 
$$
\# (\psi - \Id_A)^{-1} (a) = \# (\psi - \Id_A)^{-1} (0) = \deg (\psi- \Id_A)
$$  
and therefore, we may and do assume that $\phi$ is an isogeny.  Then $\phi^l-\Id_A$ is an isogeny so by theorem \ref{a}, 
$$
\# \Ker (\phi^l - \Id_A) = \deg ( \phi^l - \Id_A)
$$
 as $L|_{A_i}$ will be ample for every simple abelian subvariety $A_i \subset A$.  Further, $\phi|_{A_i} ^* L|_{A_i} \cong L|_{A_i} ^{\otimes q}$, since $\phi^* (L) \cong L^{\otimes q}$.  Therefore
$$
\# \{ P \in A(k) : \phi^l (P) = P \} = \deg(\phi^l - \Id_A)
$$
and the proposition follows by noticing that $\deg(\phi^l - \Id_A) = (q^l -1)^g$.
\end{proof}

We give an example showing that we can not remove the condition that $\phi^* L \cong L^{\otimes q}$ for some $q >1$.  

\begin{example}
Let $A$ be an abelian surface, $E$ and elliptic curve and let $\phi: A \times E \to A \times E$ be given by $\phi = [1] \times [4]$, where $[1] : A \to A$ and $[4] : E \to E$.  Then the degree of $\phi$ is $16$.  However, $\phi$ can not be polarized by an ample line bundle $L$.  For if this were the case, $\deg(\phi)= q^3$ for some $q \in \mathbb{N}$, which is impossible.  Furthermore, we see that $\phi$ will pointwise fix $A$.  It is possible to produce more interesting examples by taking $A$ to have complex multiplication.
\end{example}

\begin{lemma}\label{perab}
Let $A$ be an abelian variety of dimension $g$ and $L$ an ample line bundle on $A$.  Let $\phi$ be an endomorphism of $A$ with $\phi^* L \cong L^{\otimes q}$ for some $q >1$.  Suppose $Y \subset A$ is a closed normal subvariety with $\phi^m (Y) =Y$ for some $m >0$.  Then $Y$ is the translate of an abelian subvariety.
\end{lemma}

\begin{proof}
As $\phi^m |_Y : Y \to Y$ is the restriction of a positive degree endomorphism with $(\phi^m)^*L|_Y \cong L|_Y ^{\otimes q}$, for some ample line bundle $L$ on $A$, the Kodaira dimension of $Y$ can not be positive by theorem \ref{fakh}.  First we dispense with the case that $Y$ has negative Kodaira dimension.  In this case, $Y$ will be uniruled by proposition \ref{unir} and therefore, contain a rational curve.  As $Y$ is a subvariety of an abelian variety, this is impossible.  

It now suffices to consider the case when $Y$ has Kodaira dimension $0$.  Because $Y$ is a subvariety of $A$, the Albanese map $\alpha_Y : Y \to \Alb (Y)$ is injective.  By \cite[Corollary 10.6]{Ueno}, $\alpha_Y$ will also be surjective.  Let $R_{\alpha_Y}$ denote the ramification divisor of $\alpha_Y$.  Then $K_Y \sim R_{\alpha_Y} \sim 0$ as $K_{\Alb(Y)}=0$.  Appealing to lemma \ref{deblem}, we see that $Y$ is an abelian variety and hence a translate of an abelian subvariety of $A$.
\end{proof}




We are now ready to give the proof of theorem \ref{b}, which follows almost directly from the previous two lemmas. \newline

\noindent{\it Proof of theorem \ref{b} } Let $A$, $\phi$, $L$, and $Y$ be as in the statement of theorem \ref{b}.  First assume that $\phi$ is an isogeny.  By lemma \ref{perab}, $Y$ is the translate of an abelian subvariety $V$, say $Y = V +Q$, where $Q \in A(k)$.  By applying \cite[Corollary 2.2]{Fakhruddin} to a suitable power of the map $\phi^m$ along with lemma \ref{pointwise}, we see that $\phi^m$ will have a fixed point, and since $\phi^m$ is an isogeny, $Q$ must also be a fixed point of $\phi^m$.  It now suffices to consider the case where $Y$ is an abelian subvariety.  However, now we may apply proposition \ref{third} to conclude.

\section{Varieties with Kodaira Dimension $0$}\label{sec:4}

In this final section we prove theorem \ref{c}.  The outline of the proof is similar to the the proof of theorem \ref{b}.  Namely, first we prove a fixed point estimate for an endomorphism $\phi : X \to X$, where $X$ is a smooth projective variety over $\mathbb{C}$ with $\kappa (X)=0$ and $\phi^* L \cong L^{\otimes q}$ for some ample line bundle $L$ on $X$ and $q >1$.  

The idea of the proof is that by using theorem \ref{fakh}, we can restrict the geometry of $X$, so that $X \cong A / G$, where $A$ is an abelian variety and $G$ is a finite subgroup of $\Aut(G)$.  Then we may apply theorem \ref{b} to $A$ after appealing to the universal cover of $X$ to produce an endomorphism on $A$ that descends to $\phi$.  This is where we require $k=\mathbb{C}$.  We will then be able to recover a fixed point estimate for the map $\phi$ on $X$.  From the strong geometric conditions imposed on a smooth subvariety $Y \subset X$ with $\phi^m (Y) =Y$ for some $m >1$, we are then able to reduce to the case when $Y =X$, which is similar to lemma \ref{perab}.

\begin{proposition}
Let $X$ be a smooth projective variety of dimension $n$ over $\mathbb{C}$ of Kodaira dimension zero.  Let $L$ be an ample line bundle on $X$ and $\phi : X \to X$ an endomorphism with $\phi^* L \cong L^{\otimes q}$ for some $q >1$.  Then
$$
\# \{ P \in X(\mathbb{C}) : \phi^{l} (P) = P \} \geq  q^{nl} (C+o(1))
$$
as $l \to \infty$, where $C$ is a constant that depends only on $X$.
\end{proposition}

\begin{proof}
By theorem \ref{fakh}, $X$ is the quotient of an abelian variety by the action of a finite group.  In other words, $X \cong A /G$, where $A$ is an abelian variety and $G$ is a subgroup of $\Aut(A)$ with an action $G \times A \to A$.  Let $\pi : A \to X$ be the quotient map, which is finite, again by theorem \ref{fakh}.  By \cite[Proposition 2.1.1]{Zhang}, the universal cover $\widetilde{X}$ of $X$ induces a morphism $\psi : A \to A$  such that $\pi \circ \psi = \phi \circ \pi$.

Since $\pi: A \to X$ is finite, we have that $\pi^* L$ is ample on $A$.  Using that $\pi \circ \psi = \phi \circ \pi$ we then get that $\psi^* \pi^* L \cong (\pi^* L)^{\otimes q}$.  Therefore, $\psi$ has degree $q^n$.  By proposition \ref{third}, we have $\Fix (\psi) = (q -1)^n$.  Let $Q \in \Fix (\psi)$ and $g \in G$ we denote the action of $g$ on $Q$ by $g \cdot Q$.  If $G\cdot Q$ denotes the orbit of $Q$, then since $\pi(G \cdot Q)$ is a point in $X$ because $\pi \circ \psi = \phi \circ \pi$ we see that $\pi(G \cdot Q)$ is a fixed point of $\phi$.  Since each orbit of a fixed point contains at most $|G|$ fixed points, we obtain a map
$$
\Fix (\psi) \to \Fix (\phi) \textrm{ \    \ given by \    \ } G \cdot Q \mapsto \pi (G \cdot Q),
$$
which is at most $|G|$ to $1$.  Therefore, 
$$
\frac{\# \Fix (\psi)}{|G|} \leq  \# \Fix (\phi).  
$$
Since the same argument holds for any $l >>0$, we have
$$
\frac{(q^l -1)^n}{|G|} \leq \# \{ P \in X (\mathbb{C}) : \phi^l (P) =P \}.
$$ 
\end{proof}

We now give the proof of our final result, theorem \ref{c}. \newline

\noindent{\it Proof of theorem \ref{c} }
Let $X$, $L$, and $\phi$ be as in the previous proposition.  Suppose $Y \subset X$ is a smooth subvariety of dimension $r$, with $\phi^m (Y) =Y$ for some $m >0$.  In this case $L|_Y$ will be an ample line bundle on $Y$, and if $\phi^m |_Y : Y \to Y$ denotes the restriction of $\phi^m$ to $Y$, then $(\phi^m |_Y)^* L|_Y \cong (L|_Y)^{\otimes q}$.  Therefore, $\kappa(Y) \leq 0$ by theorem \ref{fakh}.  If $\kappa(Y) = - \infty$, then by proposition \ref{unir} $Y$ will be uniruled which contradicts theorem \ref{fakh} as $Y$ can not contain any rational curves.  Therefore, we may assume $\kappa(Y)=0$, then by the previous proposition applied to $Y$, $\phi^m|_Y$, and $L|_Y$ the claim follows.

\begin{remark}
One can easily give upper bounds for the number of fixed points for a subvariety $Y$ appearing in theorem \ref{c}.  By \cite[Th\'eor\`eme 1]{Serre} the absolute value of the eigenvalues of the induced map 
$$
(\phi^l)^* : H^i (X, \mathbb{C}) \to H^i (X, \mathbb{C}),
$$
are all $q^{li/2}$.  By applying the Lefschetz's fixed point theorem for $l >> 0$, see for example \cite[Chapter 3.4]{Griffiths}, we have
$$
 \# \{ P \in X (\mathbb{C}) : \phi^l (P) =P \} \leq  \sum_{i=0} ^{2n} (-1)^i \Tr((\phi^l)^* | H^i (X, \mathbb{C})) \leq q^{nl} (1 + o(1))
 $$
 which concludes the proof as $H^{2n} (X, \mathbb{C}) \cong \mathbb{C}$.  This argument does not use the fact that $Y$ is a subvariety of a variety with Kodaira dimension $0$ and in fact holds more generally \cite[Section 1.2]{Zhang}.  As was pointed out in the introduction, the main difficulty with using such an argument is that it counts the fixed points with multiplicity.  
\end{remark}

{ \it Acknowledgements }
It is a pleasure to thank my advisor M. Nakamaye for his many helpful suggestions and encouragement.  I would also like to thank A. Buium and A. Saha for many useful comments and discussions.

 

\end{document}